\documentclass[12pt]{article}
\usepackage{amsmath,amssymb,latexsym}

\usepackage[english]{babel}

\newtheorem{defn}{Definition}[section]

\newtheorem{thm}[defn]{Theorem}
\newtheorem{prop}[defn]{Proposition}
\newtheorem{cor}[defn]{Corollary}

\newcommand{\h}{{\cal H}}
\newcommand{\mn}{\mathbb N}
\def\range{{\cal R}}
\def\ru{{\cal R}(U)}
\def\h{{\cal H}}
\def\sumii{\sum_{i=1}^\infty}
\newcommand{\sumgrd[3]}{\displaystyle{\sum_{#1=#2}^#3}}
\def\bp{\noindent{\bf Proof: \ }}
\def\ep{\noindent{$\Box$}}
\newcommand{\seq[1]}{\{#1_i\}}
\newcommand{\seqgr[1]}{\{#1_i\}_{i=1}^\infty}
\def\newin {\,\kern-0.4em\in\kern-0.15em}
\def\newsubset {\kern-0.2em\subset\kern-0.2em}

\def\v{\vspace{.1in}}

\title{Perturbation of frames in Banach spaces}
\author{Diana T. Stoeva }
\begin{document}
\maketitle \pagestyle{myheadings} \markboth{Stoeva}{perturbations}

\begin{abstract}
In this paper we consider perturbation of $X_d$-Bessel sequences, $X_d$-frames, Banach frames, atomic decompositions and $X_d$-Riesz bases in separable Banach spaces. Equivalence between some perturbation conditions is investigated.
\end{abstract}

\section{Introduction}

From practical point of view, it is very important to know what happens 
with a frame for a Hilbert space $\h$, when frames' elements are changed. Throughout the years, different conditions for closeness of two frames are investigated, looking for weaker and weaker assumptions.
The first perturbation results on Hilbert frames appeared in \cite{C95pert}, where it is proved that if $\seqgr[g]$ is a frame for
$\h$ with lower bound $A$ and $\seqgr[f]$ satisfies the condition 
$\sumii\|g_i-f_i\|_\h^2<A$, 
then $\seqgr[f]$ is also a frame for $\h$. 
Recall that $\seqgr[g]\subset \h$ is called a {\it frame for the Hibert space $\h$ with bounds $A,B$} if $0<A\leq B<\infty$ and
$A \|h\|^2\leq \sumii | \langle h,g_i\rangle|^2\leq B\|h\|^2$ for every $h\in\h$.
Further perturbation results on frames  
with weaker assumptions appeared in 
[2, 5\,-\,9]. 
Perturbation of sequences, satisfying the upper frame inequality, is considered in \cite{B}. Such perturbation results are important for investigation of multipliers, which are very useful in signal processing.

As far as it is known to the author, the best condition for 
perturbation of frames up to now
is obtained by Casazza and Christensen \cite{CC}:

\begin{thm} {\rm\cite{CC}} \label{thcc}
Let $\seqgr[g]$ be a frame for $\h$ with bounds $A,B$ and let $\phi_i\newin\h$, $i\newin\mn$. If there exist constants
$\lambda_1, \lambda_2, \mu \geq 0$, such that
$\max(\lambda_1 +\frac{\mu}{\sqrt{A}}, \lambda_2)<1$
and
\begin{equation} \label{nerPW3}
\left\| \sumgrd[i]{1}{n} c_i (g_i - \phi_i ) \right\| \leq
\mu \left( \sumgrd[i]{1}{n} |c_i|^2 \right)^{\frac{1}{2}} + 
\lambda_1 \left\| \sumgrd[i]{1}{n} c_i g_i\right\| 
+ \lambda_2 \left\| \sumgrd[i]{1}{n} c_i \phi_i\right\| 
\end{equation}
for all finite scalar sequences
 $\{c_1, c_2, ..., c_n\} \ (n\in \mn)$, then $\seq[\phi]$ is a frame for $\h$ with bounds
\begin{equation} \label{pertcc} A\left(1-\frac{\lambda_1 +\lambda_2 +\mu/\sqrt{A}}{1+\lambda_2}\right)^2, \ \ \
B\left(1+\frac{\lambda_1 +\lambda_2 +\mu/\sqrt{B}}{1-\lambda_2}\right)^2. 
\end{equation}
\end{thm}

Motivated by (\ref{nerPW3}), Sun \cite{Sun} have considered perturbation of $G$-frames, which are sequences in Hilbert spaces, generalizing frames.

In the present paper we generalize condition (\ref{nerPW3}) to Banach spaces (see (\ref{bcond}) and (\ref{bcondxd})) and obtain perturbation results for some generalizations of frames to Banach spaces. The paper is organized as follows. Section \ref{nd} contains notation and needed results. Perturbation of $X_d$-Bessel sequences, $X_d$-frames, Banach frames, atomic decompositions and $X_d$-Riesz bases is the topic of Section \ref{sperturb}. For any kind of the above sequences,
we use the perturbation condition (\ref{bcond}) and determine appropriate additional assumptions on the constants $\mu, \lambda_1, \lambda_2$.
Some of the results in this section generalize results from \cite{CC,CH}.
Section \ref{sequiv} concerns connection between some conditions for closeness. 
Equivalences with simpler perturbation conditions are proved:
for $X_d$-Bessel sequences and $X_d$-frames, the $\mu$-term in (\ref{bcond}) is essential and the other two additions in (\ref{bcond}) can be omitted; 
for $X_d$-Riesz bases and Banach frames - both the $\mu$-term and the $\lambda_2$-term are essential, the $\lambda_1$-term can be omitted in some cases; for atomic decompositions - the $\mu$- and $\lambda_1$- terms can be reduced to one term.

\section{Notation, definitions and needed results}\label{nd}
Throughout the paper, $X$ and $Y$ denote Banach spaces, $X^*$ - the dual of $X$, $X_d$ - Banach space of scalar sequences. 
Recall that $X_d$ is called:
{\it $BK$-space}, if the coordinate functionals are continuous;
{\it $CB$-space}, if the canonical vectors form a Schauder basis for $X_d$;
{\it $RCB$-space}, if it is reflexive $CB$-space. 
The canonical basis of a $CB$-space is denoted by $\seqgr[e]$.

\begin{prop} \label{bkxdstar} {\rm \cite[p.\,201]{KA}}
If $X_d$ is a $CB$-space, then $X_d^\circledast \mathrel{\mathop:}=\{
\{G e_i\}_{i=1}^\infty : G\in X_d^* \}$ with the norm $\|\{G e_i\}_{i=1}^\infty\|_{X_d^\circledast}\mathrel{\mathop:}=\|G\|_{X_d^*}$ is a $BK$-space isometrically
isomorphic to $X_d^*$. 
\end{prop}

Throughout the paper, when $X_d$ is a $CB$-space, 
$X_d^*$ is identified with $X_d^\circledast$.
As usual, a scalar sequence $\seqgr[d]$ is called {\it finite}, when it has only finitely many non-zero elements. 
The notion {\it operator} is used for a linear mapping. 
It is said that an operator $F$ is defined from $X$ {\it onto} $Y$ if its range $\range({F})$ coincides with $Y$.
An operator $G$, given by $G \seqgr[c]\mathrel{\mathop:}=\sumii c_i g_i$ 
($g_i\newin Y, i\newin\mn$), is called {\it well defined from $X_d$ into $Y$} if the series $\sumii c_i g_i$ converges in $Y$ for every $\seqgr[c]\in X_d$.
The notation $\seqgr[g]\subset Y$ is used with the meaning $g_i\newin Y$, $\forall i\newin\mn$. 
If the index set of a sequence or a sum is omitted, the set $\mn$ should be understood.

Let us recall the definitions of the sequences, whose perturbations are
investigated in the present paper.

\begin{defn} \label{defxdfr}
Let $X_d$ be a $BK$-space and $\seq[g]\subset X^*$. If 
\begin{itemize}
\item[\rm{($a$)}] $ \{ g_i(f) \} \in X_d$, \ $\forall f\newin X$,
\item[\rm{($b$)}] $\exists \ \,
B\in (0,\infty)
\ \ : \ \
 \|\{g_i(f)\}\|_{X_d} \leq B \|f\|_X,  \ \forall f\newin X$,
\end{itemize}
then $\seq[g]$ is called an $X_d$-Bessel sequence for $X$ with bound $B$.

If $\seq[g]$ is an $X_d$-Bessel sequence for $X$ with bound $B$ and there exists $A\in(0,B]$ such that $A\|f\|_X\leq \|\{g_i(f)\}\|_{X_d}$ for every $f\in X$, then $\seq[g]$ is called an $X_d$-frame for $X$ with bounds $A,B$. An $\ell^p$-frame is called a $p$-frame.

When $\seq[g]$ is an $X_d$-frame for $X$ and
there exists a bounded operator $S:X_d\to X$ such that $S\{g_i(f)\}=f$, $\forall f\in X$,
then $(\seq[g], S)$ is called a {\it Banach frame for $X$ with respect to $X_d$} and $S$ is called a Banach frame operator for $\seq[g]$; $\seq[g]$ is also called a {\it Banach frame for $X$ w.r.t. $X_d$}.

When $\seq[g]$ is an $X_d$-frame for $X$ and there exists $\seq[f]\subset X$ such that $f=\sum g_i(f) f_i$, $\forall f\in X$, then $(\seq[g],\seq[f])$ is called an atomic decomposition of $X$ with respect to $X_d$.
\end{defn}

\begin{defn} \label{q} 
A sequence $\seq[g]\subset Y$ is called an {\it $X_d$-Riesz basis for $Y$ with bounds A,B}, if it is complete in $Y$, $0<A\leq B<\infty$ and
\begin{eqnarray} \label{bo}
A \|\seq[c]\|_{X_d} \le \left\| \sum c_i g_i \right\|_{Y} \le
B\|\seq[c]\|_{X_d}, \ \forall \seq[c]\in X_d. 
\end{eqnarray}
\end{defn}

Note that when $X_d$ is a $CB$-space, validity of (\ref{bo}) for all finite scalar sequences $\seqgr[c]$ implies validity of (\ref{bo}) for all $\seqgr[c]\in X_d$.

\v
While in the Hilbert space setting ($X$-Hilbert space and $X_d=\ell^2$) the concepts {\it $X_d$-frame, Banach frame} and {\it atomic decomposition} lead to a same one, namely - {\it frame for $X$}, in the Banach space setting this is not so.

1. {\it $X_d$-frame $\nRightarrow$ Banach frame w.r.t. $X_d$;

\ \ \  $X_d$-frame $\nRightarrow$ atomic decomposition w.r.t. $X_d$}:

\noindent Casazza has proved that there exist $p$-frames, which do not give rise to atomic decompositions:
{\it 
For every $p\neq 2$, $p\in(1,\infty)$, there exist a Banach space $X$ and a $p$-frame $\seq[g]\subset X^*$ for $X$ such that there is no family $\seq[f]\subset X$ satisfying $f= \sum g_i (f) f_i,  \forall f\newin X$.}
Moreover, this $p$-frame is not a Banach frame for $X$ w.r.t. $\ell^p$ (see the equivalence of {\rm(iii)} and {\rm(v)} in \cite[Proposition 3.4]{CCS}, valid for $CB$-spaces $X_d$).

2. {\it Banach frame for $X$ $\nRightarrow$ atomic decomposition of $X$}:

\noindent A sequence $\seq[g]$ is a Banach frame for $X$ if and only if $\seq[g]$ is total on $X$ i.e., if and only if $g_i(x)=0, \forall i\in \mn$, implies $x=0$ (for one of the directions see \cite[Lemma 2.6]{CCS}, the other direction is clear). Not every total sequence on $X$ give rise to atomic decomposition of $X$:  if $\seq[z]$ denotes an orthonormal basis for a Hilbert space $\h$, then the sequence $\{e_i+e_{i+1}\}$ is total and thus, it is a Banach frame for $\h$ w.r.t. appropriate $BK$-space $X_d$, however, $\{e_i+e_{i+1}\}$ does not give rise to atomic decomposition of $\h$, see \cite[Example 2.8]{CCS}.

3. {\it Atomic decomposition of $X$ w.r.t.\,$X_d$ $\nRightarrow$ Banach frame for $X$ w.r.t.\,$X_d$}:

\noindent
If $X=c_0$, $X_d={\ell}^{\infty}$ and $\seqgr[g]$ denotes the sequence of the coefficient functionals, associated to the canonical basis $\seqgr[z]$ of $c_0$, then it is clear that $(\seq[g],\seq[z])$ is an atomic decomposition of $X$ w.r.t. $X_d$. However, $\seq[g]$ is not a Banach frame for $X$ w.r.t. $X_d$ (see \cite[Example 2.3 and Proposition 3.4]{CCS}).

4. {\it  $X_d^*$-Riesz basis $\&$ $X_d$-$RCB$  $\Rightarrow$ Banach frame and atomic decomp.}:

\noindent If $X_d$ is an $RCB$-space and $\seq[g]$ is an $X_d^*$-Riesz basis for $X^*$, then $\seq[g]$ is a Banach frame for $X$ w.r.t. $X_d$  and there exists $\seq[f]$ such that $(\seq[g],\seq[f])$ is an atomic decomposition of $X$ w.r.t. $X_d$ (see \cite{Srbasis}).
Clearly, the converse does not hold in general.

\v Let $\seq[g]\subset X^*$. The operators $U$ and $T$ given by 
$$
Uf=\{g_i(f)\}, \, f\in X,
\ \mbox{and} \ T\seq[d]= \sum d_ig_i, 
$$
are called the {\it analysis operator for $\seq[g]$} and the {\it synthesis operator for $\seq[g]$}, respectively. We will use the following assertions.

\begin{prop} {\rm \cite{CCS}}  \label{prop22} Let $X_d$ be a $CB$-space (resp. $RCB$-space). 
A family $\seq[g] \subset X^* $ is an $X_d^*$-Bessel sequence (resp. $X_d$-Bessel sequence) for $X$ with bound $B$ if and only if the synthesis operator $T$
is well defined and hence bounded from  $X_d$ (resp. $X_d^*$)
  into $X^*$ and $ \|T\| \le B$.
\end{prop}

\begin{prop} {\rm \cite{Sthesis}} \label{exp}
Let $X_d$ be an $RCB$-space and $\{g_i\}\subset X^*$. The sequence $\{g_i\}$ is an $X_d$-frame for $X$ if and only if the synthesis operator $T$ is well defined and hence bounded from $X_d^*$ onto $ X^*$.
\end{prop}

For the sake of completeness we give the proof.

\bp 
Let $T$ and $U$ denote the synthesis and the analysis operator for $\seq[g]$, respectively.
 By Proposition \ref{prop22},  $\seq[g]$ is an $X_d$-Bessel sequence for $X$ if and only if $T$ is well defined from $X_d^*$ into $X^*$.

Let now $\seq[g]$ be an $X_d$-Bessel sequence for $X$.
 Consider arbitrary $F\newin X_d^*$ and the corresponding sequence $\{F(e_i)\}\newin X_d^\circledast$ (see Proposition \ref{bkxdstar}). For every $f\newin X$,
 $U^*(F) \,(f)=F(U f)=
 \sum g_i(f) \,F(e_i)=T \{F(e_i)\} \,(f)$ and hence $U^*(F)=T \{F(e_i)\} $. Thus, we can write $U^*=T$.
 Since $X_d$ is reflexive and $X$ is isomorphic to the closed
 subspace $\range(U)$ of $X_d$, $X$ is also reflexive and hence $T^*=U$.
By \cite{Heuser}, the operator $U^*$ is surjective if and only if $U$ has a bounded inverse defined on $\range(U)$, which by \cite{KA} is equivalent to the validity of the lower $X_d$-frame inequality.
\ep

\begin{prop} {\rm \cite{CCS}}  \label{propconv} Let $X_d$ be a $BK$-space and $\seq[g] \subset X^* $ be an $X_d$-frame for $X$.
If there exists $\seq[f]\subset X$ such that $\sum c_i f_i$ converges in $X$ for every $\seq[c]\in X_d$ and $f=\sum c_i f_i$, $\forall f\in X$,
then $\seq[g]$ is a Banach frame for $X$ w.r.t. $X_d$;
when $X_d$ is a $CB$-space, the converse also holds.
\end{prop}

 It is well known that a bounded operator $G:X
\rightarrow X$ on a Banach space $X$, for which $\|G-Id_X\| < 1$,
has a bounded inverse. 
An improved version of this result is given by Casazza and Christensen:

\begin{prop} {\rm\cite{CC}} \label{imprinv}
Let $G:X \rightarrow X$ be an operator. Assume that
there exist constants $\lambda_1\in [0,1), \lambda_2 \in [0,1)$ such that
\begin{equation*}
\|Gx-x\| \leq \lambda_1 \|x\| + \lambda_2 \|Gx\|, \ \forall x \in
X.
\end{equation*}
\noindent Then $G$ is bounded with bounded linear inverse $G^{-1}:X\to X$ and 
\begin{equation*}
\frac{1-\lambda_2}{1+\lambda_1} \|x\| \leq \|G^{-1}x\|
 \leq \frac{1+\lambda_2}{1-\lambda_1} \|x\|, \ \forall x\in X. 
\end{equation*}
\end{prop}

\section{Perturbation results} \label{sperturb}

Throughout the rest of the paper we assume that $\seq[g]\subset X^*$, $\seq[\phi]\subset X^*$, $\seq[f]\subset X$, $\seq[\psi]\subset X$. 

Let $X_d$ be an $RCB$-space and $\seq[g]$ be an $X_d$-Bessel sequence for $X$. 
Note that $\sum c_i g_i$ is not necessarily convergent for all $\seq[c]\in X_d$ and therefore, in general we can not generalize (\ref{nerPW3}) using $X_d$-norm of $\seq[c]$ instead of $\ell^2$-norm.
Motivated by Proposition \ref{prop22}, which implies that $\sum d_i g_i$ converges for all $\seq[d]\newin X_d^*$, we generalize condition (\ref{nerPW3}) using $X_d^*$-norm of scalar sequences. 
Thus, we consider perturbation condition in the following form:
\begin{itemize}
\item[$(\mathcal{P}^*)$] \ $\exists \ \mu\geq 0, \lambda_1\geq 0, \lambda_2\geq 0,$ such that 
\begin{equation} \label{bcond}
\left\| \sum d_i (\phi_i - g_i) \right\|_{X^*} \!\!\le
 \mu \left\| \seq[d]\right\|_{X_d^*} + \lambda_1 \left\| \sum d_i g_i \right\|_{X^*}\!\!
+ \lambda_2 \left\| \sum d_i \phi_i \right\|_{X^*}\!\!
\end{equation}
for all finite scalar sequences $\seq[d]$.
\end{itemize}

By analogue, when $\seq[f]$ is an $X_d^*$-Bessel sequence for $X^*$,  
we consider perturbation condition of the form:
\begin{itemize}
\item[$(\mathcal{P})$] \ $\exists \ \mu\geq 0, \lambda_1\geq 0, \lambda_2\geq 0,$ such that
\begin{equation} \label{bcondxd}
\left\| \sum c_i (\psi_i - f_i) \right\|_{X} \le
\mu \left\| \seq[c]\right\|_{X_d} +
\lambda_1 \left\| \sum c_i f_i \right\|_{X}
+ \lambda_2 \left\| \sum c_i \psi_i \right\|_{X}
\end{equation} 
for all finite scalar sequences $\seq[c]$.
\end{itemize}

When $\seq[g]$ is an $X_d$-Bessel sequence (resp. $\seq[f]$ is an $X_d^*$-Bessel sequence) with bound $B$ and $(\mathcal{P}^*)$ (resp.  $(\mathcal{P})$) holds, denote $$\Delta=\frac{ B(\lambda_1+\lambda_2) + \mu }{1-\lambda_2}\,.$$

\subsection{Perturbation of ${\bf X_d}$-Bessel sequences}

We begin with a result on the upper $X_d^*$-frame condition.

\begin{prop}   \label{con2xd}
Let $X_d$ be a $CB$-space and $\seq[f]$ be an $X_d^*$-Bessel sequence for $X^*$ with bound $B$. Assume that 
\begin{itemize}
\item[$\mathcal{A}_{1}:$] \ \
$(\mathcal{P})$ holds with $\lambda_2<1$.
\end{itemize}
Then $\seq[\psi]$ is an $X_d^*$-Bessel sequence for $X^*$ with bound 
$\widetilde{B}=B+\Delta$, 
 the inequality in (\ref{bcondxd}) holds for all $\seq[c] \in X_d$ and $\{\psi_i-f_i \}$ is an $X_d^*$-Bessel sequence for $X^*$ with bound $\Delta$.
\end{prop}
\bp By the triangle inequality and (\ref{bcondxd}), for every finite sequence $\seq[c]$ we have
\begin{equation*}
(1-\lambda_2) \left\|\sum c_i \psi_i \right\|_{X}
\leq      \mu \left\| \seq[c]\right\|_{X_d} +
   (\lambda_1+1)  \left\| \sum c_i f_i \right\|_{X}.
\end{equation*}
Using Proposition \ref{prop22}, for every $\seq[c]\in X_d$ and every $n> m>0$ one can conclude that
\begin{equation}\label{finew}
0\leq \left\| \displaystyle{\sum_{i=m}^n} c_i \psi_i \right\|_{X}
\leq \frac{ B(\lambda_1+1) + \mu}{1-\lambda_2} \left\| \sum_{i=m}^n c_i e_i\right\|_{X_d}\underset{n,m\to 0}{\longrightarrow} 0.
\end{equation}
Therefore, the series $\sum c_i \psi_i$ converges for every $\seq[c] \in X_d$ and
\begin{equation*}
 \left\| \displaystyle{\sum} c_i \psi_i \right\|_{X}
\leq \frac{ B(\lambda_1+1) + \mu}{1-\lambda_2} \left\| \sum c_i e_i\right\|_{X_d}, \ \, \forall \seq[c]\newin X_d.
\end{equation*}
Now Proposition \ref{prop22} implies that $\{\psi_i\}$ is an $X_d^*$-Bessel sequence for $X$ with  bound 
$\widetilde{B}=\frac{B(\lambda_1+1) + \mu }{1-\lambda_2}=B+\Delta$.

It also follows easy that $\sum c_i (\psi_i - f_i) $ converges  for all $\seq[c] \in X_d$ and the inequality in (\ref{bcondxd}) holds for all $\seq[c] \in X_d$. Moreover, 
$$
\left\| \sum c_i (\psi_i - f_i) \right\|_{X} \le
(\lambda_1 \, B + \mu + \lambda_2 \, \widetilde{B})\, \left\| \seq[c]\right\|_{X_d}
= \Delta \left\| \seq[c]\right\|_{X_d}, \forall \seq[c]\in X_d,
$$
which by Proposition \ref{prop22} implies that $\{\psi_i-f_i\}$ is an $X_d^*$-Bessel sequence for $X^*$ with bound $\Delta $.
\ep

\begin{cor}   \label{con2}
Let $X_d$ be an $RCB$-space, $\seq[g]$ be an $X_d$-Bessel sequence for $X$ with bound $B$ and 
\begin{itemize}
\item[$\mathcal{A}_{2}:$] \ \
$(\mathcal{P}^*)$ holds with $\lambda_2<1$.
\end{itemize}
Then $\seq[\phi]$ is an $X_d$-Bessel sequence for $X$ with  bound 
$\widetilde{B}=B+\Delta$, 
 the inequality in (\ref{bcond}) holds for all $\seq[d]\in X_d^*$
  and $\{\phi_i-g_i \}$ is an $X_d$-Bessel sequence for $X$ with bound $\Delta$.
\end{cor}

\subsection{Perturbation of ${\bf X_d}$-frames}

While a frame for a Hilbert space $\h$ is also a Banach frame for $\h$ and give rice to atomic decomposition of $\h$, the concepts {\it $X_d$-frame, Banach frame} and {\it atomic decomposition} are not the same in the Banach space setting. We consider perturbation of sequences of all these kinds.

\begin{prop}   \label{movepfrp5xd}
Let $X_d$ be a $CB$-space and $\seq[f]$ be an $X_d^*$-frame for $X^*$  with bounds $A, B$. 
Assume that 
\begin{itemize}
\item[$\mathcal{A}_{3}:$] \ \ 
$(\mathcal{P})$ holds with \, $\mu+\lambda_2 (A+B)+ \lambda_1 B<A$.
\end{itemize}
Then $\{ \psi _i \}$ is an $X_d^*$-frame for $X^*$  with bounds $\widetilde{A}=A- \Delta$, $\widetilde{B}=B+\Delta$.
\end{prop}

\bp 
The assumptions in $\mathcal{A}_{3}$ imply that 
$\lambda_2 <1$ and thus, 
by Proposition \ref{con2xd}, $\{\psi_i\}$ is an $X_d^*$-Bessel sequence for $X^*$ with bound $\widetilde{B}=B+\Delta$. Again by Proposition \ref{con2xd}, $\{\psi_i-f_i\}$ is an $X_d^*$-Bessel sequence for $X^*$ with bound $\Delta$
and therefore,
\begin{equation*}
\|\{g(\psi_i)\}\|_{X_d^*}\geq\|\{g(f_i)\}\|_{X_d^*} - \|\{g(\psi_i - f_i)\}  \|_{X_d^*}  \geq
(A- \Delta) \|g\|_{X^*}.
\end{equation*}
where $A- \Delta=\frac{A-\mu-\lambda_2(A+B)-\lambda_1 B}{1-\lambda_2}>0$, 
 which proves the lower $X_d^*$-frame inequality.
\ep

\begin{cor}   \label{con22}
Let $X_d$ be an $RCB$-space, 
 $\seq[g]$ be an $X_d$-frame for $X$ with bounds $A$,$B$ and 
\begin{itemize}
\item[$\mathcal{A}_{4}$\rm{:}] \ \
$(\mathcal{P^*})$ holds with $\mu+\lambda_2 (A+B)+ \lambda_1 B<A$.\end{itemize}
Then $\seq[\phi]$ is an $X_d$-frame for $X$ with bounds $\widetilde{A}=A- \Delta$, $\widetilde{B}=B+\Delta$.
\end{cor}

{\bf Note} 
 If one would like to perturb $X_d$-Bessel sequences or $X_d$-frames, keeping the new bounds close to the original ones - with difference smaller then $\varepsilon$, one can add the restriction $\Delta<\varepsilon$.

\subsection{Perturbation of Banach frames}

If $(\seq[g], S)$ is a Banach frame for $X$, there are two possibilities for perturbation:

1. perturb the operator $S$;

2. perturb the sequence $\seq[g]$.

\noindent Casazza and Christensen \cite{CC} have investigated perturbation of the Banach frame operator $S$:

\begin{thm} {\rm \cite{CC}}
Let $(\seq[g], S)$ be a Banach frame for $X$ with respect to $X_d$ with bounds $A,B$ and let the bounded operator $\widetilde{S}:X_d\to X$ satisfies the condition
\begin{itemize}
\item[$\mathcal{A}_5$\rm{:}]
 \ $\exists\, \beta_1, \beta_2, \nu \geq 0$ such that 
$\max(\beta_2, \beta_1+\nu B)<1$ and 

$\|S c-\widetilde{S} c\|_X\leq \nu \|c\|_{X_d}+
\beta_1 \|S c\|_X+
\beta_2 \|\widetilde{S} c\|_X,$ \ $\forall \, c\in X_d$.
\end{itemize}
Then there exists a sequence $\seq[\theta]\subset X^*$ such that $(\seq[\theta], \widetilde{S})$ is a Banach frame for $X$ with respect to $X_d$ with bounds $A\frac{1-\beta_2}{1+\beta_1+\nu B}$, $B\frac{1+\beta_2}{1-(\beta_1+\nu B)}$.
\end{thm}

We consider perturbation of the Banach-frame sequence.

\begin{thm}   \label{wmove}
Let $X_d$ be an $RCB$-space and $(\seq[g], S)$ be a Banach frame for $X$ with respect to $X_d$ with bounds $A,B$. Assume that 
\begin{itemize}
\item[$\mathcal{A}_{6}$\rm{:}] \ \
$(\mathcal{P^*})$ holds with 
 $\max(\lambda_2, \lambda_1+\mu \|S\|)<1$. 
\end{itemize}
 Then there exists an operator $\widetilde{S}:X_d\to X$ such that $(\seq[\phi], \widetilde{S})$ is a Banach frame for $X$ with respect to $X_d$ with bounds 
 $\widetilde{A}=\frac{1- (\mu \|S\|+\lambda_1)}{(1+\lambda_2)\|S\|}$, 
$\widetilde{B}=B+\Delta$.
 \end{thm}

\bp By Corollary \ref{con2}, $\seq[\phi]$ is an $X_d$-Bessel sequence for $X$ with bound $B+\Delta$ and hence, the operator $\widetilde{T}$, given by $\widetilde{T} \seq[d]:=\sum d_i \phi_i$, is well defined from $X_d^*$ into $X^*$ (see Proposition \ref{prop22}). 
Consider the sequence $\seq[f]:=\{S e_i\}$. By the isometrical isomorphism of $X_d^*$ and $X_d^\circledast$, $S^*(g) 
\in X_d^*$ corresponds to 
$\{g(f_i)\}=\{S^*(g)\,(e_i)\} \in X_d^\circledast, \ \forall g \in X^*$.
Moreover, $\seq[f]$ is an  $X_d^*$-frame for $X^*$ such that $g=\sum g(f_i)g_i$ for all $g \in X^*$ (see the proof of \cite[Proposition 3.4]{CCS}). 
Now we use an idea from \cite[Theorem 4]{CC}, namely, to apply Proposition \ref{imprinv} with appropriate operator $G$.
Consider the bounded operator $\widetilde{T} S^* : X^*
\rightarrow X^*$. 
Let $g \in X^*$. By Corollary \ref{con2}, the
inequality in (\ref{bcond}) holds for all $\seq[d] \in X_d^*$; applying
this inequality to 
the sequence $\{g(f_i)\}$, we get
\begin{eqnarray*}
\left\| g - \widetilde{T} S^* g \right\|_{X^*} &=&
\left\| g - \widetilde{T}  \{ g(f_i) \}  \right\|_{X^*}=
\left\| \sum g(f_i)g_i - \sum g(f_i) \phi_i \right\|_{X^*} \\
& \le &  \mu \left\|\{g(f_i)\}\right\|_{X_d^\circledast}  + \lambda_1 \left\| \sum g(f_i) g_i \right\|_{X^*} + \lambda_2 \left\| \sum g(f_i) \phi_i \right\|_{X^*}\\
&\le&  (\mu \|S\|+\lambda_1) \|g\|_{X^*} + \lambda_2 \left\| \widetilde{T} S^* g \right\|_{X^*}.
\end{eqnarray*}
By Proposition \ref{imprinv}, the operator $\widetilde{T} S^*$ 
is invertible and the inverse is bounded with
\begin{equation} \label{upperinv2}
\|(\widetilde{T} S^*)^{-1}\| \leq \frac{1+\lambda_2}{1- (\mu \|S\|+\lambda_1)}.
\end{equation}
Thus, every $g\in X^*$ can be written as $g =\widetilde{T}S^*(\widetilde{T}
S^*)^{-1} \,g$, which implies that $\widetilde{T}$ is onto $X^*$. Hence, by Proposition \ref{exp}, $\{\phi_i\}$ is an $X_d$-frame for $X$.

Since $X_d$ is reflexive and $\seq[g]$ is an $X_d$-frame for $X$, the space $X$ is isomorphic to a closed subspace of $X_d$ and thus, $X$ is also reflexive. 
Let $\widetilde{S}$ denote the bounded operator $((\widetilde{T} S^*)^{-1})^*S:X_d\to X^{**}=X$ and let $\widetilde{U}$ denote the analysis operator for $\seq[\phi]$. Note that $\widetilde{U}=\widetilde{T}^*$ (see the proof of Proposition \ref{exp}).
Therefore, $$\widetilde{S}\{\phi_i(f)\}=\widetilde{S}\widetilde{U}f=
((\widetilde{T} S^*)^{-1})^*S \widetilde{T}^*f=
f,  \, \forall f\in X,$$ and hence, $(\seq[\phi], \widetilde{S})$ is a Banach frame for $X$ with respect to $X_d$. Moreover, for every $f\in X$ we have 
\begin{equation*}
\|f\|=\|\widetilde{S} \{\phi_i(f)\}\|\leq 
\|\widetilde{S}\|\, \|\{\phi_i(f)\}\|
\leq 
\|S\|\,\frac{1+\lambda_2}{1- (\mu \|S\|+\lambda_1)} \, \|\{\phi_i(f)\}\|
\end{equation*}
and therefore, $\frac{1- (\mu \|S\|+\lambda_1)}{(1+\lambda_2)\|S\|} $ is a lower bound for $\seq[\phi]$.
\ep

\v {\bf Remark} If $\seq[g]$ is a frame for a Hilbert space $\h$ with bounds $A,B$, and $T_d$ denotes the synthesis operator for the canonical dual of $\seq[g]$, then ($\seq[g], T_d$) is a Banach frame for $\h$ with respect to $\ell^2$ with bounds $\sqrt{A},\sqrt{B}$. In this case Theorem \ref{wmove} gives Theorem \ref{thcc} - the perturbation conditions and the bounds are the same. Therefore, Theorem \ref{wmove} generalizes Theorem \ref{thcc}.

\begin{cor}   \label{wmove2}
Let $X_d$ be an $RCB$-space, $\seq[g]$ be an $X_d$-frame for $X$ with bounds $A,B$ and $P$ be a bounded projection from $X_d$ onto $\ru$, where $U$ denotes the analysis operator for $\seq[g]$.
Assume that 
\begin{itemize}
\item[$\mathcal{A}_{7}$\rm{:}] \
$(\mathcal{P^*})$ holds with 
 $\max(\lambda_2, \lambda_1+\mu \frac{\|P\|}{A})<1$.  
\end{itemize}
 Then there exists an operator $\widetilde{S}:X_d\to X$ such that $(\seq[\phi], \widetilde{S})$ is a Banach frame for $X$ with respect to $X_d$ with bounds 
$\widetilde{A}=\frac{A}{\|P\|}\frac{1- (\mu \frac{\|P\|}{A}+\lambda_1)}{(1+\lambda_2)}$,
$\widetilde{B}=B+\Delta$.
 \end{cor}
 \bp The operator $U$ has bounded inverse $U^{-1}$ with $\|U^{-1}\|\leq 1/A$.  Moreover, $S:=U^{-1}P$ is a Banach frame operator for $\seq[g]$ and $\|S\|\leq \|P\|/A$. 
 Thus, $\mathcal{A}_{6}$ holds and Theorem \ref{wmove} implies that $\seq[\phi]$ is a Banach frame for $X$ w.r.t. $X_d$ with bounds $\frac{1-(\mu \|S\|+\lambda_1)}{(1+\lambda_2)\|S\|}\geq \frac{A}{\|P\|}\frac{1-(\mu\frac{\|P\|}{A}+\lambda_1)}{1+\lambda_2}$, $B+\Delta$.
\ep

\subsection{Perturbation of ${\bf X_d}$-Riesz bases}

\begin{prop}   \label{wmoveriesz}
Let $X_d$ be an $RCB$-space and $\seq[f]$ be an 
$X_d$-Riesz basis for $X$ 
with bounds $A,B$. Assume that 
\begin{itemize}
\item[$\mathcal{A}_{8}$\rm{:}] \ $(\mathcal{P})$ holds with 
 $\max(\lambda_2, \lambda_1+\mu /A)<1$. 
 \end{itemize}
Then $\seq[\psi]$ is an $X_d$-Riesz basis for $X$ with
bounds 
$\widetilde{A}=A- \frac{A(\lambda_1+\lambda_2)+ \mu}{1+\lambda_2}$, $\widetilde{B}=B+\Delta$.
\end{prop}
\bp 
First note that $X$ is isomorphic to $X_d$ (see \cite[Proposition 3.4]{Srbasis}) and thus, $X$ is also reflexive.
 By \cite[Proposition 4.7]{Srbasis},  $\seq[f]$ is an $X_d^*$-frame for $X^*$ with bounds $A,B$ and the analysis operator $U$ for $\seq[f]$ is injective with $\range(U)=X_d^*$. 
Corollary \ref{wmove2}, applied with $P$ being the Identity operator on $X_d^*$, implies that  
$\seq[\psi]$ is an $X_d^*$-frame for $X^*$ with bounds $A\,\frac{1- (\frac{\mu}{A}+\lambda_1)}{(1+\lambda_2)}$, $B+\Delta$.
 By Proposition \ref{exp}, $\seq[\psi]$ is complete in $X^{**}=X$. By Proposition \ref{prop22}, $\{\psi_i\}$ satisfies the upper
$X_d$-Riesz basis inequality with the same upper bound $B+\Delta$.
For the lower $X_d$-Riesz basis inequality, note that
$\lambda_1<1$ 
and for every $\{c_i\} \in X_d$ one has
\begin{equation*}
\left\| \sum c_i f_i \right\| - \left\| \sum c_i \psi_i\right\| 
\le \lambda_1 \left\| \sum c_i f_i \right\| + \mu \left\|\seq[c]\right\| + \lambda_2 \left\| \sum c_i \psi_i \right\|,
\end{equation*}
\noindent which implies that
\begin{eqnarray*}
(1+\lambda_2) \left\| \sum c_i \psi_i\right\|_{X^*} &\geq&
(1- \lambda_1) \left\| \sum c_i f_i \right\|_{X^*} - \mu  \left\|\seq[c]\right\| \\
&\geq& \left( (1-\lambda_1)A - \mu \right)  \left\|\seq[c]\right\|.
\end{eqnarray*}
Therefore $\seq[\psi]$ satisfies the lower $X_d$-Riesz basis condition with bound $\widetilde{A}=\frac{(1-\lambda_1)A - \mu}{1+\lambda_2}=A- \frac{A(\lambda_1+\lambda_2)+ \mu}{1+\lambda_2}.$
\ep

\v Concerning the above proposition, note that if $A-\Delta>0$, then 
$A- \Delta$ is also a lower bound for the $X_d$-Riesz basis $\seq[\psi]$, but $A- \frac{ A(\lambda_1+\lambda_2) + \mu }{1+\lambda_2}$ is closer to the optimal one.

{\bf Remark} If $\seq[g]$ is a Riesz bazis for a Hilbert space $\h$, then Proposition \ref{wmoveriesz} becomes \cite[Corollary 2]{CC}. Therefore, 
Proposition \ref{wmoveriesz} is a generalization of \cite[Corollary 2]{CC} to Banach spaces.

\subsection{Perturbation of atomic decompositions}

Perturbation of atomic decompositions under ($\mathcal{P}$) with $\lambda_2=0$ is considered in \cite{CH}. Below we add the $\lambda_2$-term and 
assume (\ref{bcondxd}) only for a subspace of $X_d$.
The reason to work with a subspace of $X_d$ is the following.
Let $X_d$ be a $BK$-space and
$(\seq[g], \seq[f])$ be an atomic decomposition of $X$ with respect to $X_d$. Clearly, $\sum c_i f_i$ converges in $X$ for every $\seq[c]=\{g_i(f)\}, f\in X$. However,  
$\sum c_i f_i$ does not need to converge for all $\seq[c]$ in $X_d$, $\seq[f]$ might be not an $X_d^*$-Bessel sequence - an example can be found in \cite{LO}.
That is why below we only assume that (\ref{bcondxd}) holds for $\seq[c]=\{g_i(f)\}, f\in X$, not necessarily for all $\seq[c]\in X_d$.
Note that convergence of $\sum c_i f_i$ in $X$ for all $\seq[c]\in X_d$ implies that $\seq[g]$ is a Banach frame for $X$ w.r.t. $X_d$  (see Proposition \ref{propconv}).
 
\begin{prop}   \label{atdmove2}
Let $X_d$ be an $RCB$-space and $(\seq[g], \seq[f])$ be an atomic decomposition of $X$ with respect to $X_d$ with bounds $A,B$.
Assume that 
\begin{itemize}  
\item[$\mathcal{A}_{9}$\rm{:}] \ $\exists$ $\mu\geq 0, \lambda_1\geq 0, \lambda_2\geq 0,$ $\max(\lambda_2, \lambda_1+\mu B)<1$, such that (\ref{bcondxd}) holds for every $\seq[c]=\{g_i(f)\}, f\in X$.
\end{itemize}
 Then there exists a sequence $\seq[\theta]\subset X^*$ such that $(\seq[\theta], \seq[\psi])$ is an atomic decomposition of $X$ with respect to $X_d$ with bounds 
 $\widetilde{A}=A\frac{1-\lambda_2}{1+\lambda_1+\mu B}$, 
$\widetilde{B}=B\frac{1+\lambda_2}{1-(\lambda_1+\mu B)}$.
 \end{prop}

\bp By $\mathcal{A}_{9}$, $\sum g_i(f)\psi_i$ converges in $X$ for every $f\in X$. Thus, we can consider the operator $G:X\to X$ given by $Gf:=\sum g_i(f) \psi_i$, $f\in X$. 
For every $f\in X$, one has
\begin{eqnarray*}
\left\| f-Gf \right\|_{X} &=& 
\left\| \sum g_i(f) f_i - \sum g_i(f) \psi_i \right\|_{X} 
\\
&\le&
\mu \left\| \{g_i(f)\} \right\|_{X_d} +
\lambda_1 \left\| \sum g_i(f) f_i \right\|_{X}
+ \lambda_2 \left\| \sum g_i(f) \psi_i \right\|_{X}\\
&\le& (\lambda_1+B\mu) \left\| f \right\|_{X}
+ \lambda_2 \left\| Gf \right\|_{X}.
\end{eqnarray*}
By Proposition \ref{imprinv}, $G$ is bounded with bounded inverse.
For $i\in\mn$, define $\theta_i:= (G^{-1})^*g_i\in X^*$. 
For every $f\in X$, one has 
$$\{\theta_i(f)\}=\{g_i(G^{-1}f)\}\in X_d,$$
\begin{equation*}\|\{\theta_i(f)\}\|_{X_d}=\|\{g_i(G^{-1}f)\}\|_{X_d}
\left\{ 
\begin{array}{l}
\leq B\|G^{-1}f\|\leq B\,\frac{1+\lambda_2}{1-(\lambda_1+\mu B)}\,\|f\|\\
\geq A\|G^{-1}f\|\geq A\,\frac{1-\lambda_2}{1+\lambda_1+\mu B}\,\|f\|
\end{array}
\right.,
\end{equation*}
$$f=G(G^{-1}f)=\sum g_i(G^{-1}f)\psi_i =\sum \theta_i(f)\psi_i.$$
This concludes the proof.
\ep

\begin{thm}  \label{atdmove3}
Let $X_d$ be an $RCB$-space and $(\seq[g], \seq[f])$ be an atomic decomposition of $X$ with respect to $X_d$ with bounds $A,B$. Assume that 
\begin{itemize}
\item[$\mathcal{A}_{10}$\rm{:}] \ $\exists$ $\mu\geq 0, \lambda_1\geq 0, \lambda_2\geq 0,$ $\max(\lambda_2, \lambda_1+\mu B)<1$, such that (\ref{bcondxd}) holds for every $\{c_i\}_{i=1}^n=\{g_i(f)\}_{i=1}^n, n\in\mn, f\in X$.
\end{itemize}
 Then the conclusion of Proposition \ref{atdmove3} holds.
 \end{thm}
\bp
Let $f\in X$. First prove the convergence of $\sum g_i(f)\psi_i$ in $X$. 
By $\mathcal{A}_{10}$, for $m>n$ we obtain
\begin{equation*}\label{psiconv2}
(1-\lambda_2) \left\| \sum_{i=n+1}^m g_i(f) \psi_i \right\|_{X} \le
\mu \left\| \sum_{i=n+1}^m g_i(f)e_i\right\|_{X_d} +
(\lambda_1 +1)\left\| \sum_{i=n+1}^m g_i(f) f_i \right\|_{X} \to 0
\end{equation*}
when $m,n\to\infty$, 
which implies that $\sum g_i(f)\psi_i$ converges in $X$. 
Furthermore, consider 
$$ \left\| \sum_{i=1}^n g_i(f) (\psi_i - f_i) \right\|_{X} \!\!\!\le
\mu \left\| \sum_{i=1}^n g_i(f) e_i\right\|_{X_d} \!\!\!\!+
\lambda_1 \left\| \sum_{i=1}^n g_i(f) f_i \right\|_{X}
\!\!+ \lambda_2 \left\| \sum_{i=1}^n g_i(f) \psi_i \right\|_{X}
$$
and take limit as $n\to\infty$. This implies that (\ref{bcondxd}) holds
for every $\seq[c]=\{g_i(f)\}$, $f\in X$, and therefore, $\mathcal{A}_{9}$ holds.
The rest follows from Proposition \ref{atdmove2}.
\ep

\v 
Concerning the above theorem, note that if in addition $\sum c_i f_i$ converges in $X$ for all $\seq[c]\in X_d$
and (\ref{bcondxd}) holds for all finite sequences $\seq[c]$, then $\sum c_i \psi_i$ converges in $X$ for all $\seq[c]\in X_d$. Thus, the following assertion follows easy from Theorem \ref{atdmove3} and Proposition \ref{propconv}.

\begin{cor} 
Let the assumptions of Theorem \ref{atdmove3} hold and let  in $\mathcal{A}_{10}$ it is assumed that (\ref{bcondxd}) holds for all finite sequences $\seq[c]$.
If $\seq[g]$ is a Banach frame for $X$ w.r.t. $X_d$, then there exists a sequence $\seq[\theta]\subset X^*$, which is a Banach frame for $X$ w.r.t. $X_d$, and such that $(\seq[\theta], \seq[\psi])$ is an atomic decomposition of $X$ with respect to $X_d$.
\end{cor}

\section{Perturbation conditions} \label{sequiv}

In this section we investigate how essential are the terms 
in (\ref{bcond}) and (\ref{bcondxd}). 
We prove that some of the terms can be omitted and simpler perturbation conditions can be used.

\subsection{$X_d^*$-Bessel sequences}

We begin with observation that the terms with $\lambda_1$ and $\lambda_2$ in (\ref{bcondxd}) are not essential for perturbation of $X_d^*$-Bessel sequences. 

\begin{prop} \label{bequiv} Let $X_d$ be a $CB$-space and 
 $\seq[f]$ be an $X_d^*$-Bessel sequence for $X^*$ with  bound
$B$. Then $\mathcal{A}_{1}$ is equivalent to the following conditions for closeness:
\begin{itemize}
\item[$\mathcal{A}_{11}${\rm :}] \
$\exists \, \widetilde{\mu} \geq 0$ such that  
\begin{equation}\label{bs}
\left\| \sum c_i (\psi_i - f_i) \right\|_{X} \le
\widetilde{\mu} \| \{c_i\}\|_{X_d}
\end{equation}
 for all finite sequences $\seq[c]$ (and hence for all $\seq[c]\in X_d$).
 \item[$\mathcal{A}_{12}${\rm :}] $\{g(\psi_i - f_i)\} \in X_d^*$, $\forall g\in X^*$, and 
$\exists \, \widetilde{\mu} \geq 0$ such that  
\begin{equation}\label{b2b}
\| \{g(\psi_i - f_i)\} \|_{X_d^*} \le
\widetilde{\mu} \| g\|_{X^*}, \ \forall g\in X^*.
\end{equation}
\item[$\mathcal{A}_{13}${\rm :}] $\{g(\psi_i)\} \in X_d^*$, $\forall g\in X^*$, and 
$\exists \, \mu\geq 0, \lambda_1\geq 0, \lambda_2\in [0,1),$  such that  
\begin{equation}\label{bcond6cb}
\| \{g(\psi_i - f_i)\} \|_{X_d^*} \le
\mu \| g\|_{X^*} + \lambda_1 \| \{ g(f_i) \|_{X_d^*} + \lambda_2 \| \{g( \psi_i)\} \|_{X_d^*}, \ \forall g\in X^*.
\end{equation}
\end{itemize}
\end{prop}

\bp 

$\mathcal{A}_{11} \Leftrightarrow \mathcal{A}_{12} $: If one of the conditions $\mathcal{A}_{11} $ and $\mathcal{A}_{12}$ holds with $\widetilde{\mu}=0$, then $\seq[f]\equiv\seq[\psi]$ and thus the other condition also holds with $\widetilde{\mu}=0$. For the cases when $\widetilde{\mu}\neq 0$, the equivalence of $\mathcal{A}_{11} $ and $\mathcal{A}_{12}$ follows from Proposition \ref{prop22}.

$\mathcal{A}_{11}\Rightarrow \mathcal{A}_{1}$ and $\mathcal{A}_{12} \Rightarrow \mathcal{A}_{13}$: obvious.

$\mathcal{A}_{1} \Rightarrow \mathcal{A}_{11}$: Let $\mathcal{A}_{1}$ hold.
By Proposition \ref{con2}, $\{\psi_i - f_i\}$ is an $X_d^*$-Bessel sequence for $X^*$ with bound $\Delta$ and now Proposition \ref{prop22} implies that (\ref{bs}) holds with $\widetilde{\mu}=\Delta\geq 0$.

$\mathcal{A}_{13}  \Rightarrow \mathcal{A}_{12}$: Let $\mathcal{A}_{13}$ hold. For every $g\in X^*$, it follows that $\{g(\psi_i - f_i)\} \in X_d^*$ and
$$
\| \{g(\psi_i)\} \|_{X_d^*} \leq 
(\lambda_1+1 ) \| \{g(f_i)\} \|_{X_d^*} + \mu \|g\|_{X^*} + \lambda_2 \| \{g(\psi_i)\} \|_{X_d^*}.
$$
Hence,
$$(1-\lambda_2) \| \{g(\psi_i)\} \|_{X_d^*} \leq 
((\lambda_1+1 )B + \mu)\|g\|_{X^*}, 
$$ which implies that
$\seq[\psi]$ is an $X_d^*$-Bessel sequence for $X^*$ with  bound $\widetilde{B}=\frac{(\lambda_1+1 )B + \mu}{1-\lambda_2}=B+\Delta$. Therefore,
$$\| \{g(\psi_i - f_i)\} \|_{X_d^*} \leq 
(\mu + \lambda_1 B +\lambda_2 \widetilde{B})\|g\|_{X^*}=\Delta \|g\|_{X^*},
$$
i.e.  $\mathcal{A}_{12}$ holds.
\ep

\v Note that if $\seq[f]$ and $\seq[\psi]$ are $X_d^*$-Bessel sequences for $X^*$, then $\{f_i-\psi_i\}$ is also an $X_d^*$-Bessel sequence for $X^*$ and thus, $\mathcal{A}_{12}$ is a natural perturbation assumption for $X_d^*$-Bessel sequences.

\subsection{$X_d^*$-frames}

As in the $X_d^*$-Bessel case, the terms with 
$\lambda_1$ and $\lambda_2$ are not essential for $X_d^*$-frames. Note that going from the $X_d^*$-Bessel case to the $X_d^*$-frame case, we add the restriction $\widetilde{\mu}<A$ to $\mathcal{A}_{11}$ and $\mathcal{A}_{12}$. This restriction is essential - if $\mathcal{A}_{11}$ holds with $\widetilde{\mu}=A$, then $\seq[\psi]$ is not needed to be an $X_d^*$-frame for $X$.

\begin{prop} \label{xdfrequiv}
Let $X_d$ be a $CB$-space and $\seq[f]$ be an $X_d^*$-frame for $X^*$ with bounds $A,B$. Then $\mathcal{A}_{3}$ is equivalent to the following conditions for closeness:
\begin{itemize}
\item[$(a)$\rm{:}] \ $\mathcal{A}_{11}$ holds with $\widetilde{\mu}<A$.  
\item[$(b)$\rm{:}] \ $\mathcal{A}_{12}$ holds with $\widetilde{\mu}<A$.
\item[$(c)$\rm{:}] \ $\exists \, \mu\geq 0, \lambda_1\geq 0, \lambda_2>0, \mu + \lambda_2(A+B) + \lambda_1 B<A$, such that  (\ref{bcond6cb}) holds.
\end{itemize}
\end{prop}
\bp 

$(a)\Leftrightarrow (b)$: Follows in the same way as in Proposition \ref{bequiv}.

$\mathcal{A}_{3} \Rightarrow (a)$: Let $\mathcal{A}_{3}$ hold and hence, $\Delta\in [0,A)$. 
By Proposition \ref{con2}, $\{\psi_i - f_i\}$ is an $X_d^*$-Bessel sequence for $X^*$ with bound $\Delta$ and now Proposition \ref{prop22} implies that (\ref{bs}) holds with $\widetilde{\mu}=\Delta$.

$(a) \Rightarrow \mathcal{A}_{3}$ and $(b) \Rightarrow (c)$: obvious.

$(c) \Rightarrow (b)$: Let $(c)$ hold and hence, $\Delta\in [0,A)$. By the proof of $\mathcal{A}_{13}\Rightarrow \mathcal{A}_{12}$, it follows that (\ref{b2b}) holds with $\widetilde{\mu}=\Delta$.
\ep

\v In a similar way as above, assertions concerning equivalence of perturbation conditions for $X_d$-Bessel sequences and $X_d$-frames 
can be written.

\subsection{Banach frames and ${\bf X_d}$-Riesz bases}

For perturbation of $X_d$-Riesz bases and Banach frames,
the $\lambda_1$-term can be omitted in certain cases. 
For general Banach frames it is clear that $\|S\|\geq 1/B$. The following proposition concerns the case 
when the equality holds.

\begin{prop} Let $X_d$ be an $RCB$-space and $(\seq[g], S)$ be a Banach frame for $X$ with respect to $X_d$ with bounds $A,B$. If $\|S\|=1/B$, then  
$\mathcal{A}_{6}$ is equivalent to the following condition:
 \begin{itemize}
\item[$\widetilde{\mathcal{A}_{6}}$\rm{:}] \
$\exists \widetilde{\mu}\geq 0, \widetilde{\lambda}_2\geq 0, \max(\widetilde{\lambda}_2, \widetilde{\mu}/ B)<1$, such that 
\begin{equation*}\label{rcond}
\left\| \sum d_i (\phi_i - g_i) \right\|_{X} \le
 \widetilde{\mu} \| \{d_i\}\|_{X_d^*}
+ \widetilde{\lambda_2} \left\| \sum d_i \phi_i \right\|_{X}
\end{equation*}
for all finite sequences $\{d_i\}$ (and hence for all $\seq[d]\in X_d^*$).
\end{itemize}
\end{prop}
\bp
For the implication ($\mathcal{A}_{6}$ $\Rightarrow$ $\widetilde{\mathcal{A}_6}$), take $\widetilde{\mu}:=\mu+\lambda_1 B$, $\widetilde{\lambda}_2:=\lambda_2$. 
The other implication is obvious.
\ep

\begin{prop} \label{r1}
Let $X_d$ be an $RCB$-space, $\seq[f]$ be an $X_d$-Riesz basis for $X$ with bounds $A,B$. 
If $A=B$, then  
$\mathcal{A}_{8}$ is equivalent to 
the following condition:
 \begin{itemize}
\item[$\widetilde{\mathcal{A}_{8}}$\rm{:}] \
$\exists \, \widetilde{\mu}\geq 0, \widetilde{\lambda}_2\geq 0, \max(\widetilde{\lambda}_2, \widetilde{\mu} /A)<1$, such that 
\begin{equation*}
\left\| \sum c_i (\psi_i - f_i) \right\|_{X} \le
 \widetilde{\mu} \| \{c_i\}\|_{X_d}
+ \widetilde{\lambda_2} \left\| \sum c_i \psi_i \right\|_{X}
\end{equation*}
for all finite sequences $\{c_i\}$ (and hence for all $\seq[c]\in X_d$).
\end{itemize}
\end{prop}
\bp
For the implication ($\mathcal{A}_{8}$ $\Rightarrow$ $\widetilde{\mathcal{A}_8}$), take $\widetilde{\mu}:=\mu+\lambda_1 B$, $\widetilde{\lambda}_2:=\lambda_2$. The other implication is obvious.
\ep

\subsection{Atomic decompositions} 

Concerning perturbation of atomic decompositions, recall that $\mathcal{A}_{9}$ requires validity of (\ref{bcondxd}) only for $\seq[c]=\{g_i(f)\}$. In this case the $\lambda_1$-term is $\lambda_1\|f\|$. Thus, we can replace the $\lambda_1$- and the $\mu$- terms by one term with $\|f\|$:

\begin{prop} \label{a1} 
Let $X_d$ be an $RCB$-space and 
$(\seq[g],\seq[f])$ be an atomic decomposition of $X$ with respect to $X_d$ with upper bound $B$. 
Then $\mathcal{A}_{9}$ is equivalent to the following condition:
 \begin{itemize} 
\item[$\widetilde{\mathcal{A}_{9}}$\rm{:}] \
$\exists \, \widetilde{\mu}\geq 0, \widetilde{\lambda}_2\geq 0, \max(\widetilde{\lambda}_1, \widetilde{\lambda}_2)<1$, such that 
\begin{equation} \label{adeq}
\left\| \sum g_i(f) (\psi_i - f_i) \right\|_{X} \le
 \widetilde{\lambda}_1 \| f\|_{X}
+ \widetilde{\lambda_2} \left\| \sum g_i(f) \psi_i \right\|_{X}
\end{equation}
for every $f\in X$.
\end{itemize}
\end{prop}

\bp
For the implication ($\mathcal{A}_{9}$ $\Rightarrow$ $\widetilde{\mathcal{A}_9}$), take $\widetilde{\lambda}_1=\lambda_1 +\mu B$, $\widetilde{\lambda}_2=\lambda_2$. 
The other implication is clear.
\ep

\vspace{.1 in}
\noindent Diana T. Stoeva\\
Department of Mathematics \\
University of Architecture, Civil Engineering and Geodesy \\
Blvd. Christo Smirnenski 1,\\
1046 Sofia \\
Bulgaria \\
Email: stoeva\_\,fte@uacg.bg


\begin{thebibliography}{150}

\frenchspacing

\bibitem{B} Balazs, P.:  Basic definition and properties of Bessel multipliers.
{\it J. Math. Anal. Appl.} {\bf 325}, No. 1 (2007), 571-585. 


\bibitem{CC} Casazza, P. and Christensen, O.:
 Perturbation of operators and applications to frame theory.
{\it J. Fourier Anal. Appl.} {\bf 3} (1997), No.5, 543-557. 



\bibitem{CCS}  Casazza, P.G., Christensen, O. and Stoeva, D.T.:
 Frame expansions in separable Banach spaces.  {\it J. Math. Anal.
Appl.}, \textbf{307} (2005), 710--723.

\bibitem{C95pert} Christensen, O.:
 Frame perturbations.
{\it Proc. Am. Math. Soc.} {\bf 123}, No.4 (1995), 1217-1220.

\bibitem{C95} Christensen, O.:
 A Paley-Wiener theorem for frames. {\it Proc. Am. Math. Soc.} {\bf
123}, No.7 (1995), 2199-2202.


\bibitem{C96pert}
Christensen, O.: 
Moment problems and stability results for frames with applications to irregular sampling and Gabor frames.
{\it Appl. Comput. Harmon. Anal.} {\bf 3}, No.1 (1996), 82-86.




\bibitem{Cpertbook}
Christensen, O.: {\it Perturbation of frames and applications to Gabor frames.} In \lq\lq Gabor analysis and algorithms. Theory and applications.\rq\rq \
Eds. H. G. Feichtinger and T. Strohmer, Birkh$\ddot{a}$user, 1998, 193-209. 

\bibitem{CH} Christensen, O. and Heil, C.:
 Perturbations of Banach frames and atomic decompositions.
{\it Math. Nachr.}, {\bf 185} (1997), 33-47.

\bibitem{FZ} Favier, S.J. and Zalik, R.A.: On the stability of frames and Riesz bases. {\it Appl. Comput. Harmon. Anal.} {\bf 2}, No.2 (1995), 160-173.

\bibitem{Heuser}
H. Heuser, {\small\it Functional analysis,} (John Wiley \& Sons, 1982).

\bibitem{KA}   Kantorovich, L.\,V. and   Akilov, G.\,P.
{\it Functional analysis in normed spaces}, Pergamon press, Oxford, 1964.

\bibitem{LO} S. Li  and H. Ogawa, Pseudo-duals of frames with applications, {\it Appl. Comput. Harm. Anal.}, {\bf 11}, No. 2 (2001), 289-304.

\bibitem{Sthesis} Stoeva, D.T.: {\it Frames and bases in Banach spaces}, PhD thesis, UACG and UCTM, 2005.

\bibitem{Srbasis} Stoeva, D.T.: $X_d$-Riesz bases in separable Banach spaces, \lq\lq Collection of papers, ded. to the 60th Anniv. of M. Konstantinov\rq\rq, BAS Publ. House, Sofia, 2008.

\bibitem{Sun} Sun, W.: Stability of $G$-frames. 
{\it J. Math. Anal. Appl.}, {\bf 326}, No. 2, p. 858-868.


\end{thebibliography}
\end{document}